\documentclass{svmult}
\usepackage{makeidx,graphicx,multicol,amssymb,amsmath}
\usepackage[bottom]{footmisc}

\newcommand{\R}{{\mathbb R}}
\renewcommand{\S}{{\mathbb S}^{d-1}}
\newcommand{\N}{{\mathbb N}}

\newcommand{\ird}[1]{\int_{\R^d}{#1}\;dx}
\newcommand{\nrm}[2]{\|{#1}\|_{\L^{#2}(\R^d)}}
\newcommand{\nrmcnd}[2]{\|{#1}\|_{\L^{#2}(\mathcal C)}}

\newcommand{\ic}[1]{\int_{\mathcal C}{#1}\;dy}

\newcommand{\be}[1]{\begin{equation}\label{#1}}
\newcommand{\ee}{\end{equation}}
\renewcommand{\(}{\left(}
\renewcommand{\)}{\right)}
\renewcommand{\S}{{\mathbb S^{d-1}}}
\newcommand{\C}[1]{\mathsf C_{\rm #1}}
\renewcommand{\H}{\mathrm H}
\renewcommand{\L}{\mathrm L}

\newcommand{\CKN}{{(CKN)}}
\newcommand{\WLH}{{(WLH)}}

\begin{document}

\title*{About existence, symmetry and symmetry breaking for extremal functions of some interpolation functional inequalities}
\titlerunning{Existence and symmetry of extremals in some interpolation inequalities}

\author{Jean Dolbeault\and Maria J. Esteban}
\institute{Ceremade (UMR CNRS no.~7534), Universit\'e Paris-Dauphine, Place de Lattre de Tassigny, F-75775 Paris C\'edex 16, France\\
\texttt{dolbeaul@ceremade.dauphine.fr}, \texttt{esteban@ceremade.dauphine.fr}
}
\maketitle

\begin{abstract}
This article is devoted to a review of some recent results on existence, symmetry and symmetry breaking of optimal functions for Caffarelli-Kohn-Nirenberg \CKN\ and weighted logarithmic Hardy \WLH\ inequalities. These results have been obtained in a series of papers \cite{DET,0902,DDFT,DE2010,1007} in collaboration with M.~del Pino, S.~Filippas, M.~Loss, G.~Tarantello and A.~Tertikas. Here we put the highlights on a symmetry breaking result: extremals of some inequalities are not radially symmetric in regions where the symmetric extremals are linearly stable. Special attention is paid to the study of the critical cases for \CKN\ and \WLH.
\end{abstract}


\section{Two families of interpolation inequalities}\label{Sec:Interpolation}

For any dimension $d\in\N^*\,$ and any $\,\theta\in[0,1]$, let us consider the set $\mathcal D$ of all smooth functions which are compactly supported in $\R^d\setminus\{0\}$. Define the numbers
\begin{multline*}
\vartheta(p,d):=\frac{d\,(p-2)}{2\,p}\,,\;\,a_c:=\frac{d-2}2\,,\;\Lambda(a):=(a-a_c)^2\\\mbox{and}\quad p(a,b):=\frac{2\,d}{d-2+2\,(b-a)}\,.
\end{multline*}
We shall also set $\,2^*:=\frac{2\,d}{d-2}\,$ if $\,d\ge 3\,$ and $\,2^*:=\infty\,$ if $\,d=1\,$ or~Ê$\,2$. For any $\,a<a_c$, we consider the following two families of interpolation inequalities, which have been introduced in \cite{Caffarelli-Kohn-Nirenberg-84,DDFT}:
\begin{description}
\item[(CKN)] {\it Caffarelli-Kohn-Nirenberg inequalities} -- Let $b\in(a+1/2,a+1]$ and $\theta\in(1/2,1]$ if $d=1$, $b\in(a,a+1]$ if $d=2$ and $b\in[a,a+1]$ if $d\ge3$. Assume that $p=p(a,b)$, and $\theta\in[\vartheta(p,d),1]$ if $d\ge2$. Then, there exists a finite positive constant $\C{CKN}(\theta,p,a)$ such that, for any $u\in\mathcal D$,
\[
\nrm{|x|^{-b}\,u}p^2\le\C{CKN}(\theta,p,a)\,\nrm{|x|^{-a}\,\nabla u}2^{2\,\theta}\,\nrm{|x|^{-(a+1)}\,u}2^{2\,(1-\theta)}\,.
\]
\item[(WLH)] {\it Weighted logarithmic Hardy inequalities} --
Let $\gamma \ge d/4$ and $\gamma>1/2$ if $d=2$. There exists a positive constant $\C{WLH}(\gamma,a)$ such that, for any $u\in\mathcal D$, normalized by $\nrm{|x|^{-(a+1)}\,u}2=1$,
\[
\ird{\frac{|u|^2\,\log\(|x|^{2\,(a_c-a)}\,|u|^2\)}{|x|^{2\,(a+1)}}}\leq2\,\gamma\,\log\!\left[\C{WLH}(\gamma,a)\,\nrm{\,|x|^{-a}\,\nabla u}2^2\right]\,.
\]
\end{description}
\WLH\ appears as a limiting case of \CKN\ in the limit \hbox{$\theta=\gamma\,(p-2)$}, \hbox{$p\to 2_+$}. See \cite{DDFT,DE2010} for details. By a standard completion argument, these inequalities can be extended to the set
\[
\mathcal D_a^{1,2}(\R^d):=\{u\in\L^1_{\rm loc}(\R^d)\,:\,|x|^{-a}\,\nabla u\in\L^2(\R^d)\;\mbox{and}\;|x|^{-(a+1)}\, u\in\L^2(\R^d)\}\,.
\]
In the sequel, we shall assume that all constants in the inequalities are taken with their optimal values. For brevity, we shall call {\it extremals} the functions which realize equality in \CKN\ or in \WLH.

Let $\C{CKN}^*(\theta,p,a)$ and $\C{WLH}^*(\gamma,a)$ denote the optimal constants when admissible functions are restricted to the set of radially symmetric functions. {\it Radial extremals}, that is, the extremals of the above inequalities when restricted to the set of radially symmetric functions, are explicit and the values of the constants, $\C{CKN}^*(\theta,p,a)$ and $\C{WLH}^*(\gamma,a)$, are known. According to \cite{DDFT}, we have:
\begin{description}
\item[(CKN$^*$)]{\it Radial Caffarelli-Kohn-Nirenberg inequalities}:
\[\label{Eqn:OptC*}
\textstyle\frac{\C{CKN}^*(\theta,p,a)}{|\S|^{-\frac{p-2}p}}=\left[\frac{(a-a_c)^2\,(p-2)^2}{2+(2\theta-1)\,p}\right]^\frac{p-2}{2\,p}\!\left[\frac{2+(2\theta-1)\,p}{2\,p\,\theta\,(a-a_c)^2}\right]^\theta \!\left[\frac 4{p+2}\right]^\frac{6-p}{2\,p}\!\left[\frac{\Gamma\left(\frac{2}{p-2}+\frac 12\right)}{\sqrt\pi\;\Gamma\left(\frac{2}{p-2}\right)}\right]^\frac{p-2}{p}\!,
\]
and if $\theta>\vartheta(p,1)$, the best constant is achieved by an optimal radial function $u$ such that $u(r)=r^{a-a_c}\,\overline{w}(-\log r)$, where $\overline{w}$ is unique up to multiplication by constants and translations in $s$, and given by
\[
\overline{w}(s)=\big(\cosh(\lambda\,s)\big)^{-\frac2{p-2}},\quad\mbox{with}\quad{\textstyle\lambda =\frac 12\,(p-2)\,\left[\frac{(a-a_c)^2 (p+2)}{2+(2\theta-1)\,p}\right]^{\frac12}}\,.
\]
\item[(WLH$^*$)] {\it Radial weighted logarithmic Hardy inequalities}:
\end{description}
\[
\textstyle\C{WLH}^*(\gamma,a) =\left\{\begin{array}{ll}\frac{1}{4\,\gamma}\, \frac{ \left[ \Gamma \(\frac{d}{2}\) \right]^{\frac{1}{2\,\gamma}}}{(2\,\pi^{d+1}\,e)^{\frac{1}{4\,\gamma}}}\(\frac{4\,\gamma -1}{(a-a_c)^2}\)^{ \frac{4\,\gamma -1}{4\,\gamma}}\quad&\mbox{if}\;\gamma>\frac 14\,,\\
\quad\C{WLH}^*= \frac{ \left[ \Gamma \(\frac{d}{2}\) \right]^2}{2\,\pi^{d+1}\,e}\quad&\mbox{if}\;\gamma=\frac 14\,,
\end{array}\right.
\]
and if $\gamma>\frac 14$, equality in the weighted logarithmic Hardy inequality is achieved by an optimal radial function $u$ such that $u(r)=r^{a-a_c}\,w(-\log r)$, where
\[
w(s)= \frac{\tilde w(s)}{ \ic{\tilde w^2}}\quad \hbox{\rm and}\quad \tilde w(s) = \exp \(-\frac{(a-a_c)^2\,s^2}{(4\,\gamma-1)} \)\,.
\]
Moreover we have
\be{Ineq:CompRad}\begin{array}{l}
\C{CKN}(\theta,p,a)\ge\C{CKN}^*(\theta,p,a)=\C{CKN}^*(\theta,p, a_c-1)\,\Lambda(a)^{\frac{p-2}{2p}-\theta}\,,\\[6pt]
\C{WLH}(\gamma,a)\ge\C{WLH}^*(\gamma,a)=\C{WLH}^*(\gamma,a_c-1)\,\Lambda(a)^{-1+\frac 1{4\,\gamma}}\,,
\end{array}\ee
where in both cases, the inequalities follow from the definitions. 
Radial symmetry for the extremals of \CKN\ and \WLH\ implies that $\C{CKN}(\theta,p,a)=\C{CKN}^*(\theta,p,a)$ and $\,\C{WLH}(\gamma,a)=\C{WLH}^*(\gamma,a)$, while {\it symmetry breaking} means that inequalities in~\eqref{Ineq:CompRad} are strict. As we shall see later, there are cases where $\C{CKN}(\theta,p,a)=\C{CKN}^*(\theta,p,a)$ and for which radial and non radial extremal functions coexist. This may happen only for the limiting value of $a$ beyond which the equality does not hold anymore. On the contrary, when $\C{CKN}(\theta,p,a)>\C{CKN}^*(\theta,p,a)$, none of the extremals of \CKN\ is radially symmetric.

\medskip Section~\ref{Sec:Extremals} is devoted to the attainability of the best constants in the above inequalities. In Section~\ref{Sec:Symmetry} we describe we describe the best available symmetry breaking results. In Section~\ref{Sec4} we give some plots and also prove some new asymptotic results in the limit $p\to 2_+$.

\section{Existence of extremals}\label{Sec:Extremals}

In this section, we describe the set of parameters for which the inequalities are achieved. The following result is taken from~\cite{DE2010}.
\begin{theorem}[Existence based on {\it a priori} estimates]\label{Thm:Existence} Equality in \CKN\ is attained for any $p\in(2,2^*)$ and $\theta\in(\vartheta(p,d),1)$ or for $\,\theta=\vartheta(p,d)$ and $a\in(a_\star,a_c)$, for some $a_\star<a_c$. It is not attained if $p=2$, or $a<0$, $p=2^*$, $\theta=1$ and $d\ge 3$, or $d=1$ and $\theta=\vartheta(p,1)$.

Equality in \WLH\ is attained if $\,\gamma\ge 1/4\,$ and $\,d=1$, or $\,\gamma>1/2\,$ if $\,d=2$, or for $\,d\geq 3\,$ and either $\,\gamma>d/4\,$ or $\,\gamma=d/4\,$ and $\,a\in(a_{\star\star},a_c)$, for some $\,a_{\star\star}<a_c$. 
\end{theorem}
A complete proof of these results is given in \cite{DE2010}. In the sequel, we shall only give some indications on how they are established. 

\medskip First of all, it is very convenient to reformulate \CKN\ and \WLH\ inequalities in cylindrical variables. By means of the Emden-Fowler transformation
\[
s=\log|x|\in\R\,,\quad\omega=x/|x|\in\S\,,\quad y=(s,\omega)\,,\quad v(y)=|x|^{a_c-a}\,u(x)\,,
\]
inequality \CKN\ for $u$ is equivalent to a Gagliardo-Nirenberg-Sobolev inequality on the cylinder $\mathcal C:=\R\times\S$ for $v$, namely
\[
\nrmcnd vp^2\leq\C{CKN}(\theta,p,a)\(\;\nrmcnd{\nabla v}2^2+\Lambda\,\nrmcnd v2^2\)^\theta\,\nrmcnd v2^{2\,(1-\theta)}\quad\forall\;v\in\H^1(\mathcal C)
\]
with $\Lambda=\Lambda(a)$. Similarly, with $w(y)=|x|^{a_c-a}\,u(x)$, inequality \WLH\ is equivalent to
\[\label{Ineq:GLogHardy-w}
\ic{|w|^2\,\log |w|^2}\leq 2\,\gamma\,\log\left[\C{WLH}(\gamma,a)\left(\nrmcnd{\nabla w}2^2+\Lambda\right)\right]\,,
\]
for any $w\in\H^1(\mathcal C)$ such that $\nrmcnd w2=1$. {\it Notice that radial symmetry for $u$ means that $v$ and $w$ depend only on $s$.} For brevity, we shall call them {\it $s$-symmetric functions}.

\medskip On $\H^1(\mathcal C)$, consider the functional
\[\label{en1}
\mathcal E_{\theta,\Lambda}^p[v]:=\(\nrmcnd{\nabla v}2^2\!+\Lambda\,\nrmcnd v2^2\)^\theta\nrmcnd v2^{2\,(1-\theta)}\,.
\]
Assume that $d\ge3$, let $t:=\nrmcnd{\nabla v}2^2/\nrmcnd v2^2$ and $\Lambda=\Lambda(a)$. If $v$ is a minimizer of $\mathcal E_{\theta,\Lambda}^p[v]$ such that $\nrmcnd vp=1$, then, as in \cite{DE2010}, we have
\begin{multline}\label{est1}
(t+\Lambda)^\theta=\mathcal E_{\theta,\Lambda}^p[v]\,\frac{\nrmcnd vp^2}{\nrmcnd v2^2}=\frac{\nrmcnd vp^2}{\C{CKN}(\theta,p,a)\,\nrmcnd v2^2}\\
\le\frac{\mathsf S_d^{\vartheta(p,d)}}{\C{CKN}(\theta,p,a)}\,\(t\!+a_c^2\)^{\vartheta(p,d)}
\end{multline}
where $\mathsf S_d=\C{CKN}(1,2^*,0)$ is the optimal Sobolev constant, while we know from~\eqref{Ineq:CompRad}, that $\lim_{a\to a_c}\C{CKN}(\theta,p,a)=\infty$ if $d\ge 2$. This provides a bound on~$t$ if $\,\theta>\vartheta(p,d)$. 

Consider now a sequence $(v_n)_n$ of functions in $\H^1(\mathcal C)$, which minimizes $\mathcal E_{\theta,\Lambda}^p[v]$ under the constraint $\nrmcnd vp=1$. Assume therefore that $\nrmcnd{v_n}p=1$ for any $n\in\N$. If $\,\theta>\vartheta(p,d)$, an estimate similar to~\eqref{est1} asymptotically holds for $(v_n)_n$, thus providing bounds on $t_n:=\|\nabla v_n\|_{\L^2({\mathcal C})}/\| v_n\|_{\L^2({\mathcal C})}\,$ and $\,\| v_n\|_{\H^1({\mathcal C})}$, for $n$ large enough.

Then, standard tools of the concentration-compactness method allow to conclude that $(v_n)_n$ is relatively compact and converges up to translations and the extraction of a subsequence towards a minimizer of $\mathcal E_{\theta,\Lambda}^p$. The only specific idea concerning the use of concentration-compactness in this context relies on the use of the following inequality: for any $x$, $y>0$ and any $\eta\in(0,1)$,
\[\label{alg1}
(1+x)^\eta\,(1+y)^{1-\eta}\ge 1+x^\eta\,y^{1-\eta}, \;\mbox{with strict inequality unless $x=y$}\,.
\]

A similar approach holds for \CKN\ if $d=2$.

\medskip In the case of \WLH, for $\gamma>d/4\,$, the method of proof is similar to that of \CKN. The energy functional to be considered is now
\[\label{en2}
\mathcal F_\gamma[w]:=\(\nrmcnd{\nabla w}2^2\!+\Lambda\)\,\exp\left[-\frac1{2\,\gamma}\ic{|w|^2\,\log|w|^2}\right]\,.
\]
If $w$ is a minimizer of $\mathcal F_\gamma[w]$ under the constraint $\nrmcnd{w}2=1$, then we have
\[\label{est2}
\frac{t+\Lambda}{\big[\C{CKN}(1,p,\alpha)\,(t+\Lambda(\alpha))\big]^{\frac 1{2\,\gamma}\frac p{p-2}}}\le\mathcal F_\gamma[w]=\frac1{\C{WLH}(\gamma,a)}\le\frac{\Lambda(a)^{1-\frac 1{4\,\gamma}}}{\C{WLH}^*(\gamma,a_c-1)}
\]
for an arbitrary $\alpha<a_c$. The concentration-compactness method applies using the following inequality: for any $x$, $y>0$ and $\eta\in (0,1)$,
\[
\eta\,x^{1/\eta}+(1-\eta)\,y^{1/(1-\eta)}\ge x\,y\,,
\]
with strict inequality unless $x=y$ and $\eta=1/2$.

\medskip Let us now consider the {\it critical case} $\,\theta=\vartheta(p,d)$ for \CKN. Estimate~\eqref{est2} still provides {\it a priori} bounds for minimizing sequences whenever $\,a\in (a_1, a_c)$ where $a_1$ can be obtained as follows. When $\,\theta=\vartheta(p,d)$, we can rewrite~\eqref{est1} as
\[
(t+\Lambda) \le\,K\,\(t\!+a_c^2\)\quad\mbox{where}\quad\mathsf K=\frac{\mathsf S_d^\theta}{\C{CKN}(\theta,p,a)}\,.
\]
Hence we can deduce that $\, 0\leq t\le\frac{\mathsf K\,a_c^2-\Lambda}{1-\mathsf K}$, if $\mathsf K<1$ and $\,\Lambda\leq \mathsf K\,a_c^2$. These two inequalities define the constant
\be{Lambda1}
\Lambda_1:=\min\left\{\left(\frac{\C{CKN}^*(\theta,p,a_c-1)^{1/\theta}}{\mathsf S_d}\right)^\frac d{d-1},\left(\frac{a_c^2\,\mathsf S_d}{\C{CKN}^*(\theta,p,a_c-1)^{1/\theta}}\right)^d\right\}\,,
\ee
so that $t$ is bounded if $\,a\in (a_1, a_c)$ with $a_1:=a_c-\sqrt{\Lambda_1}$. See \cite{DE2010} for more details.

\medskip Such an estimate is not anymore available in the {\it critical case} for \WLH, that is, if $\gamma=d/4$, $d\ge3$. We may indeed notice that $p\le 2^*$ and $\gamma=d/4$ mean $1-\frac 1{2\,\gamma}\frac p{p-2}\le 0$. A more detailed analysis of the possible losses of compactness is therefore necessary. This can actually be done in the two critical cases, $\,\theta=\vartheta(p,d)$ for \CKN\ and $\gamma=d/4$, $d\ge3$, for \WLH. 

\medskip Let $\C{GN}(p)$ be the optimal constant in the Gagliardo-Nirenberg-Sobolev interpolation inequalities
\[
\nrm up^2\le\C{GN}(p)\,\nrm{\nabla u}2^{2\,\vartheta(p,d)}\,\nrm u2^{2\,(1-\vartheta(p,d))}\quad\forall\;u\in\H^1(\R^d)
\]
with $p\in(2,2^*)$ if $d=2$ or $p\in(2,2^*]$ if $d\ge3$. Also consider Gross' logarithmic Sobolev inequality in Weissler's form (see \cite{MR479373})
\be{LSI}
\ird{|u|^2\,\log |u|^2}\le\frac d2\,\log\(\C{LS}\,\nrm{\nabla u}2^2\)
\ee
for any $u\in\H^1(\R^d)$ such that $\nrm u2=1$, with optimal constant \hbox{$\C{LS}\!:=\frac 2{\pi\,d\,e}$}.

The gaussian function $u(x)=\frac 1{(2\,\pi)^{d/4}}\,e^{-|x|^2/4}$ is an extremal for~\eqref {LSI}. By taking $u_n(x):=u(x+n\,\mathsf e)$ for some $\mathsf e\in\S$ and any $n\in\N$ as test functions for \WLH, and letting $n\to+\infty$, we find that
\[
\C{LS}\le\C{WLH}(d/4,a)\,.
\]
If equality holds, this is a mechanism of loss of compactness for minimizing sequences. On the opposite, if $\C{LS}<\C{WLH}(d/4,a)$, we can establish a compactness result (see Theorem~\ref {MainThm2} below) which proves that, for some $a_{\star\star}<a_c$, equality is attained in \WLH\ in the critical case $\gamma=d/4$ for any $a\in(a_{\star\star},a_c)$. Indeed, we know that $\lim_{a\to a_c}\C{WLH}(d/4,a)=\lim_{a\to a_c}\C{WLH}^*(d/4,a)=\infty$.

\medskip A similar analysis for \CKN\ shows that
\[
\C{GN}(p)\le\C{CKN}(\theta,p,a)
\]
in the critical case $\theta=\vartheta(p,d)$. Exactly as for \WLH, we also have an existence result, which has been established in \cite{DE2010}, if $\C{GN}(p)<\C{CKN}(\theta,p,a)$.
\begin{theorem}[Existence in the critical cases]\label{MainThm2} With the above notations,
\begin{itemize}
\item[{\rm (i)}] if $\theta=\vartheta(p,d)$ and $\C{GN}(p)<\C{CKN}(\theta,p,a)$, then \CKN\ admits an extremal function in $\mathcal D^{1,2}_{a}(\R^d)$,
\item[{\rm (ii)}] if $\gamma=d/4$, $d\ge 3$, and $\C{LS}<\C{WLH}(\gamma,a)$, then \WLH\ admits an extremal function in $\mathcal D^{1,2}_{a}(\R^d)$. Additionnally, if $a\in(a_{\star\star}^{\scriptscriptstyle\rm WLH},a_c)$, then $\C{LS}<\C{WLH}(d/4,a)$ where $a_{\star\star}^{\scriptscriptstyle\rm WLH}$ is defined by
\[
a_{\star\star}^{\scriptscriptstyle\rm WLH}:=a_c-\sqrt{\Lambda_{\star\star}^{\scriptscriptstyle\rm WLH}}\quad\mbox{and}\quad \Lambda_{\star\star}^{\scriptscriptstyle\rm WLH}:=(d-1)\,e\,\left[\frac{\Gamma(\frac d2)^2}{2^{d+1}\,\pi}\right]^\frac1{d-1}\,.
\]
\end{itemize}\end{theorem}
The values of $\C{GN}(p)$ and $\C{CKN}(\vartheta(p,d),p,a)$ are not explicitly known if $d\ge2$, so we cannot get an explicit interval of existence in terms of $a$ for \CKN. The strict inequality of Theorem~\ref{MainThm2} (i) holds if $\C{GN}(p)<\C{CKN}^*(\vartheta(p,d),p,a)$ since we know that $\C{CKN}^*(\vartheta(p,d),p,a)\le\C{CKN}(\vartheta(p,d),p,a)$. The condition $\C{GN}(p)=\C{CKN}^*(\vartheta(p,d),p,a)$ defines a number $a_\star^{\scriptscriptstyle\rm CKN}$ for which existence is granted if $a\in(a_\star^{\scriptscriptstyle\rm CKN},a_c)$, hence proving that $a_\star\le a_\star^{\scriptscriptstyle\rm CKN}$ (if we consider the lowest possible value of $a_\star$ in Theorem~\ref{Thm:Existence}). Still we do not know the explicit value of $\C{GN}(p)$, but, since the computation of $a_1$ only involves the optimal constants among radial functions, at least we know that $a_\star^{\scriptscriptstyle\rm CKN}\leq a_1 <a_c$.

On the opposite, we know the explicit values of $\C{LS}$ and $\C{WLH}^*(d/4)$, so that the computation of the value of $a_\star^{\scriptscriptstyle\rm WLH}$, which is determined by the condition $\C{LS}=\C{WLH}^*(d/4)$, is tedious but explicit.

\medskip We may observe from the expression of \CKN\ and \WLH\ when they are written on the cylinder (after the Emden-Fowler transformation) that $\C{CKN}$ and $\C{WLH}$ are monotone non-decreasing functions of $a$ in $(-\infty,a_c)$, and actually increasing if there is an extremal. As long as it is finite, the optimal function $a_\star$ in Theorem~\ref{Thm:Existence} is continuous as a function of $p$, as a consequence of Theorem~\ref{MainThm2} and of the compactness of minimizing sequences. So, finally, $a_\star$ and $a_{\star\star}$ can be chosen such that $a_\star\le a_\star^{\scriptscriptstyle\rm CKN}$ and $a_{\star\star}\le a_{\star\star}^{\scriptscriptstyle\rm WLH}$. It is not difficult to observe that $a_\star^{\scriptscriptstyle\rm CKN}$ can be seen as a continuous, but not explicit, function of $p$ and we shall see later (in Corollary~\ref{cor1}, below) that $\lim_{p\to 2_+}a_\star^{\scriptscriptstyle\rm CKN}(p)=a_{\star\star}^{\scriptscriptstyle\rm WLH}$.

\medskip Next, note that if $\C{CKN}=\C{CKN}^*$ is known, then there are radially symmetric extremals, whose existence has been established in \cite{DDFT}. Anticipating on the results of the next section, we can state the following result which arises as a consequence of the Schwarz symmetrization method (see Theorem~\ref{Thm:Symmetry}, below, and \cite{1007}).
\begin{proposition}[Existence of radial extremals]\label{Prop} Let $d\ge3$. Then \CKN\ with $\theta=\vartheta(p,d)$ admits a radial extremal if $a\in[a_0,a_c)$ where $a_0:=a_c-\sqrt{\Lambda_0}$ and $\Lambda=\Lambda_0$ is defined by the condition
\[
\Lambda^{(d-1)/d}=\vartheta(p,d)\,\C{CKN}^*(\theta,p,a_c-1)^{1/\vartheta(p,d)}/\,\mathsf S_d\,.
\] \end{proposition}
A similar estimate also holds if $\theta>\vartheta(p,d)$, with less explicit computations. See \cite{1007} for details.

\medskip
The proof of this symmetry result follows from a not straightforward use of the Schwarz symmetrization. If
$u(x)=|x|^a\,v(x)$, \CKN\ is equivalent~to
\[
\nrm{|x|^{a-b}\,v}p^2\le\C{CKN}(\theta,p,\Lambda)\(\mathcal A-\lambda\,\mathcal B\)^\theta\,\mathcal B^{1-\theta}
\]
with {$\mathcal A:=\nrm{\nabla v}2^2$, $\mathcal B:=\nrm{|x|^{-1}\,v}2^2\;$} and $\;\lambda:=a\,(2\,a_c-a)$.

We observe that the function $B\mapsto h(\mathcal B):=\(\mathcal A-\lambda\,\mathcal B\)^\theta\,\mathcal B^{1-\theta}$ satisfies
\[
\frac{h'(\mathcal B)}{h(\mathcal B)}=\frac{1-\theta}{\mathcal B}-\frac{\lambda\,\theta}{\mathcal A-\lambda\,\mathcal B}\,.
\]
By Hardy's inequality ($d\ge3$), we know that
\[
\mathcal A-\lambda\,\mathcal B\ge\inf_{a>0}\big(\mathcal A-a\,(2\,a_c-a)\,\mathcal B\big)=\mathcal A-a_c^2\,\mathcal B>0\,,
\]
and so $h'(\mathcal B)\le 0\,$ if $\,(1-\theta)\,\mathcal A<\lambda\,\mathcal B$, which is equivalent to ${\mathcal A/\mathcal B}<\lambda/(1-\theta)$.
By interpolation ${\mathcal A/\mathcal B}$ is small if $a_c-a>0$ is small enough, for $\theta>\vartheta(p,d)$ and $d\ge3$. The precise estimate of when $\,{\mathcal A/\mathcal B}\,$ is smaller than $\,\lambda/(1-\theta)\,$ provides us with the definition of $\,a_0$.

\section{Symmetry and symmetry breaking}\label{Sec:Symmetry}

Define
\[
\underline a(\theta,p):= a_c-\frac{2\,\sqrt{d-1}}{p+2}\,\sqrt{\frac{2\,p\,\theta}{p-2}-1}\,,\quad \tilde a(\gamma):=a_c-\frac 12\sqrt{(d-1)(4\,\gamma-1)}\,,
\]
\[
a_{\rm SB}:=a_c-\sqrt{\Lambda_{\rm SB}(\gamma)}\,,\;\Lambda_{\rm SB}(\gamma):=\frac{4\,\gamma-1}8\,e\!\(\tfrac{\pi^{4\,\gamma-d-1}}{16}\)^\frac 1{4\,\gamma-1}\!\(\tfrac d\gamma\)^\frac{4\,\gamma}{4\,\gamma-1}\!\Gamma\(\tfrac d2\)^\frac 2{4\,\gamma-1}
\]
and take into account the definitions of $\,a_\star^{\scriptscriptstyle\rm CKN}\,$ and $\,a_{\star\star}^{\scriptscriptstyle\rm WLH}\,$ previously given. Thus we have the following result, which has been established in \cite{DDFT,1007}.
\begin{theorem}\label{Thm:SymmetryBreaking} Let $d\ge 2$ and $p\in(2,2^*)$. Symmetry breaking holds in \CKN\ if either $a<\underline a(\theta,p)$ and $\theta\in[\vartheta(p,d),1]$, or $a<a_\star^{\scriptscriptstyle\rm CKN}$ and $\theta=\vartheta(p,d)$.

Assume that $\gamma>1/2$ if $d=2$ and $\gamma\ge d/4$ if $d\ge3$. Symmetry breaking holds in \WLH\ if $a<\max\{\tilde a(\gamma),a_{\rm SB}\}$.\end{theorem}
When $\gamma=d/4$, $d\ge 3$, we observe that $\Lambda_{\star\star}^{\scriptscriptstyle\rm WLH}=\Lambda_{\rm SB}(d/4)<\Lambda(\tilde a(d/4))$ with the notations of Theorem~\ref {Thm:Existence} and there is symmetry breaking if $a\in(-\infty,a_{\star\star}^{\scriptscriptstyle\rm WLH})$, in the sense that $\C{WLH}(d/4,a)>\C{WLH}^*(d/4,a)$ in that interval, although we do not know if extremals for \WLH\ exist when $\gamma=d/4$ and $a<a_{\star\star}^{\scriptscriptstyle\rm WLH}.$

Concerning \CKN\ with $\,\theta\geq \vartheta(p,d)$, results of symmetry breaking for $a<\underline a(\theta,p)$ have been established first in \cite{Catrina-Wang-01,Felli-Schneider-03,DET} when $\theta=1$ and later extended in \cite{DDFT} to $\theta<1$. The main idea in case of \CKN\ is to consider the quadratic form associated to the second variation of $\mathcal E_{\theta,\Lambda}^p$, restricted to $\{v\in\H^1(\mathcal C)\,:\,\nrmcnd vp=1\}$, around a minimizer among functions depending on $s$ only and observe that the linear operator $\mathcal L_{\theta,\Lambda}^p$ associated to the quadratic form has a negative eigenvalue if $a<\underline a$. 

Because of the homogeneity in \CKN, if $v$ is a $s$-symmetric extremal, then $\lambda\,v$ is also a $s$-symmetric extremal for any $\lambda\in\R$ and $v$ is therefore in the kernel of $\mathcal L_{\theta,\Lambda}^p$. When $v$ generates $\mathrm{Ker}(\mathcal L_{\theta,\Lambda}^p)$ and all non-zero eigenvalues are positive, that is for $a\in(\underline a(\theta,p),a_c)$, we shall say that $v$ is {\it linearly} stable, without further precision. In such a case, the operator $\mathcal L_{\theta,\Lambda}^p$ has the property of {\it spectral gap}.

Results in \cite{DDFT} for \WLH, $a<\tilde a(\gamma)$, are based on the same method.

\medskip For any $a<a_\star^{\scriptscriptstyle\rm CKN}$, we have
\[
\C{CKN}^*(\vartheta(p,d),p,a)<\C{GN}(p)\le\C{CKN}(\vartheta(p,d),p,a)\,,
\]
which proves symmetry breaking. Using well-chosen test functions, it has been proved \cite{1007} that $\underline a(\vartheta(p,d),p)<a_\star^{\scriptscriptstyle\rm CKN}$ for $p-2>0$, small enough, thus also proving symmetry breaking for $a-\underline a(\vartheta(p,d),p)>0$, small, and \hbox{$\theta-\vartheta(p,d)>0$}, small. {\it This shows that in some cases, symmetry can be broken even in regions where the radial extremals are linearly stable. In Section 4, we give a more quantitative result about this (see Corollary~\ref{cor1}).}

\medskip Next we will describe how the set of parameters involved in our inequalities is cut into two subsets, both of them simply connected. They are separated by a continuous surface which isolates the symmetry region from the region of symmetry breaking. See \cite{0902,1007} for detailed statements and proofs.
\begin{theorem}\label{Thm:Symmetry} For all $d\geq 2$, there exists a continuous function $a^*$ defined on the set $\{(\theta,p)\in(0,1]\times(2,2^*)\,:\,\theta>\vartheta(p,d)\}$ such that \hbox{$\lim_{p\to 2_+}a^*(\theta,p)=-\infty$} with the property that \CKN\ has only radially symmetric extremals if $(a,p)\in(a^*(\theta,p),a_c)\times(2,2^*)$, and none of the extremals is radially symmetric if $(a,p)\in(-\infty,a^*(\theta,p))\times(2,2^*)$.

Similarly, for all $d\geq 2$, there exists a continuous function $a^{**}:(d/4,\infty)\to(-\infty,a_c)$ such that, for any $\gamma>d/4$ and $a\in [a^{**}(\gamma), a_c)$, there is a radially symmetric extremal for \WLH, while for $a<a^{**}(\gamma)$ no extremal is radially symmetric.
\end{theorem}
We sketch below the main steps of the proof. First note that as previously explained (see \cite{1007} for details), the Schwarz symmetrization allows to characterize a nonempty subdomain of $(0,a_c)\times(0,1)\ni(a,\theta)$ in which symmetry holds for extremals of \CKN, when $d\ge 3$. If $\theta=\vartheta(p,d)$ and $p>2$, there are radially symmetric extremals if $a\in[a_0,a_c)$ where $a_0$ is given in Proposition~\ref{Prop}.

Symmetry also holds if $a_c-a$ is small enough, for \CKN\ as well as for \WLH, or when $p\to 2_+$ in \CKN, for any $d\ge 2$, as a consequence of the existence of the spectral gap of $\mathcal L_{\theta,\Lambda}^p$ when $a>\underline a(\theta,p)$.

According to \cite{0902,1007}, for given $\theta$ and $p$, there is a unique $a^*\in(-\infty,a_c)$ for which there is symmetry breaking in $(-\infty,a^*)$ and for which all extremals are radially symmetric when $a\in(a^*,a_c)$. This follows from the observation that, if $v_\sigma(s,\omega):=v(\sigma\,s,\omega)$ for $\sigma>0$, then the quantity
\[
(\mathcal E_{\theta,\sigma^2\Lambda}^p[v_\sigma])^{1/\theta}-\sigma^{(2\,\theta-1+2/p)/\theta^2}\,(\mathcal E_{\theta,\Lambda}^p[v])^{1/\theta}
\]
is equal to $0$ if $v$ depends only on~$s$, while it has the sign of $\sigma-1$ otherwise. The method also applies to \WLH\ and gives a similar result for $a^{**}$.

\medskip From Theorem~\ref{Thm:SymmetryBreaking}, we can infer that radial and non-radial extremals for \CKN\ with $\theta>\vartheta(p,d)$ coexist on the threshold, in some cases.

\section{Numerical computations and asymptotic results for \CKN}\label{Sec4}

\medskip In the critical case for \CKN, that is for $\theta=\vartheta(p,d)$, numerical results illustrating our results on existence and on symmetry {\it versus} symmetry breaking have been collected in Figs.~1~and~2 below.

\subsection{Existence for \CKN}

\begin{figure}[!h]
\begin{center}
\includegraphics[height=5cm]{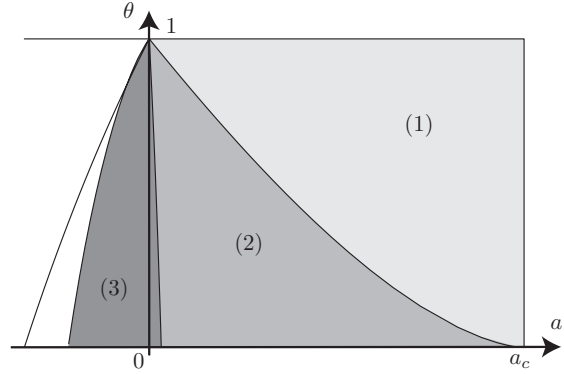}
\caption{Existence in the critical case for \CKN. Here we assume that $d=5$.}
\label{fig1}
\end{center}
\end{figure}

\noindent In Fig.~1, the zones in which existence is known are:
\begin{description}
\item (1) $a\ge a_0$: extremals are achieved among radial functions, by the Schwarz symmetrization method (Proposition~\ref{Prop}),
\item (1)+(2) $a>a_1$: this follows from the explicit {\it a priori} estimates (Theorem~\ref{Thm:Existence}); see~\eqref{Lambda1} for the definition of $\Lambda_1=(a_c-a_1)^2$,
\item (1)+(2)+(3) $a>a_\star^{\scriptscriptstyle\rm CKN}$: this follows by comparison of the optimal constant for \CKN\ with the optimal constant in the corresponding Gagliardo-Nirenberg-Sobolev inequality (Theorem~\ref{MainThm2}).
\end{description}

\subsection{Symmetry and symmetry breaking for \CKN}

\noindent In Fig.~2, the zone of symmetry breaking contains:
\begin{description}
\item (1) $a<\underline a(\theta,p)$: by linearization around radial extremals (Theorem~\ref{Thm:SymmetryBreaking}),
\item (1)+(2) $a< a_\star^{\scriptscriptstyle\rm CKN}$: by comparison with the Gagliardo-Nirenberg-Sobolev inequality (Theorem~\ref{Thm:SymmetryBreaking}).
\end{description}
In (3) it is not known whether symmetry holds or if there is symmetry breaking, while in (4), that is, for $a_0\le a<a_c$, according to Proposition~\ref{Prop}, symmetry holds by the Schwarz symmetrization.

\begin{figure}[!ht]
\begin{center}
\includegraphics[height=5cm]{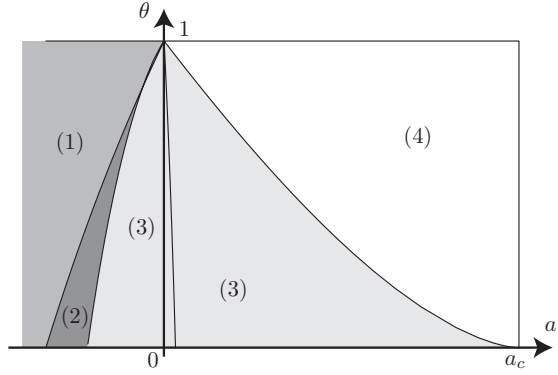}
\caption{Symmetry and symmetry breaking results in the critical case for \CKN. Here we assume that $d=5$.}
\label{fig2}
\end{center}
\end{figure}

\subsection{When \CKN\ approaches \WLH}

In the critical case $\theta=\vartheta(p,d)=d\,(p-2)/(2\,p)$, when $p$ approaches $2_+$, it is possible to obtain detailed results for \CKN\ and to compare \CKN\ and \WLH, or at least get explicit results for the various curves of Figs.~1and~2.

\smallskip\noindent 1) {\it Cases covered by the Schwarz symmetrization method.} With $\Lambda_0$ defined by $\Lambda_0^{(d-1)/d}=\vartheta(p,d)\,\C{CKN}^*(\theta,p,a_c-1)^{1/\vartheta(p,d)}/\,\mathsf S_d$, since
\[
\lim_{p\to2_+}\C{CKN}^*(\theta,p,a_c-1)^{1/\vartheta(p,d)}=\frac{(d-1)^\frac{d-1}d}{d\,(2\,e)^{1/d}\,\pi^\frac{d+1}d}\,\Gamma\big(\tfrac d2\big)^{2/d}
\]
it follows that $a_0$ defined in Proposition~\ref{Prop} by $a_0=a_c-\sqrt{\Lambda_0}$ converges to $a_c$ as $p\to2_+$.

\smallskip\noindent 2) {\it Existence range obtained by a priori estimates.} The expression of $a_1=a_c-\sqrt{\Lambda_1}$ is explicit for any $p$ and $p\mapsto\Lambda_1(p)$ has a limit $\Lambda_1(2)$ as $p\to 2_+$, which is given by
\[
\textstyle\min\left\{\frac 14\left[\frac 2e\,(d-2)^d\, (d-1)^{d-3}\(\frac {\Gamma\big(\tfrac{d}{2}\big)}{\Gamma \big(\tfrac{d-1}{2}\big)}\)^{\!2}\,\right]^{\frac 1{d-1}}\kern -12pt,\;\frac e8\,\frac{(d-2)^d}{(d-1)^{d-3}}\(\frac{\Gamma \big(\tfrac{d-1}{2}\big)}{\Gamma\big(\tfrac{d}{2}\big)}\)^{\!2}\,\right\}\,.
\]
A careful investigation shows that $\Lambda_1(2)$ is given by the first term in the above $\min$. As a function of $d$, $a_c-\sqrt{\Lambda_1(2)}$ is monotone decreasing in $(3,\infty)$ and converges to $0_+$ as $d\to\infty$.  Moreover, for all $d\ge 2$, $\Lambda_1(2)\leq \Lambda_{\star\star}^{\scriptscriptstyle\rm WLH}$, since both estimates are done among radial functions and the latter is optimal among those.

\smallskip\noindent 3) {\it Symmetry breaking range obtained by linearization around radial extremals.} Computations are explicit and it has already been observed in \cite{DDFT} that $\underline a(\theta,p)$ (see Theorem~\ref{Thm:SymmetryBreaking}) is such that $\lim_{p\to2_+}\underline a(\vartheta(p,d),p)=-1/2$.

\smallskip\noindent 4) {\it Existence range obtained by comparison with Gagliardo-Nirenberg-Sobolev inequalities.} Although the value of $\C{GN}(p)$ is not known explicitly, we can get an estimate by using a Gaussian as a test function. This estimate turns out to be sharp as $p$ approaches $2_+$. More precisely, we get a lower bound for $\C{GN}(p)$ by computing
\[
Q(p):=\frac{\nrm {u_2}p^2}{\nrm{\nabla {u_2}}2^{2\,\vartheta(p,d)}\,\nrm {u_2}2^{2\,(1-\vartheta(p,d))}}
\]
with $u_2(x):=\pi^{-d/4}\,e^{-|x|^2/2}$, which is such that
\[
\textstyle{\displaystyle \lim_{p\to2_+}}\frac{Q(p)-1}{p-2}=\frac d4\,\log\C{LS}=\frac d4\,\log\(\frac 2{\pi\,d\,e}\)\le{\displaystyle \lim_{p\to2_+}}\frac{\C{GN}(p)-1}{p-2}\,.
\]
This estimate is not only a lower bound for the limit, but gives its exact value, as shown by the following new result.
\begin{proposition}\label{Lem:Limit} With the above notations, we have
\[
{\displaystyle \lim_{p\to2_+}}\frac{\C{GN}(p)-1}{p-2}=\frac d4\,\log\C{LS}\,.
\]
\end{proposition}
Hence, in the regime $p\to2_+$, the condition which defines $a=a_\star^{\scriptscriptstyle\rm CKN}$, namely the equality $\,\C{GN}(p)=\C{CKN}^*(\vartheta(p,d),p,a)$ leads to
\begin{multline*}
1+\frac d4\,\log\C{LS}\,(p-2)+o(p-2)=\C{GN}(p)\\
=\C{CKN}^*(\vartheta(p,d),p,a)=1+\frac d4\,\log\C{WLH}^*(d/4,a)\,(p-2)+o(p-2)
\end{multline*}
(for the second line in the inequality, see \cite[Lemma 4]{DDFT}), which asymptotically amounts to solve
\[
\C{WLH}^*(d/4,a)=\C{LS}\,.
\]
In other words, we have
\[
\lim_{p\to2_+}a_\star^{\scriptscriptstyle\rm CKN}(p)=a_{\star\star}^{\scriptscriptstyle\rm WLH}\,.
\]
As a consequence, we have the following symmetry breaking result, which allows to refine an earlier result of \cite{1007} in the subcritical case and is new in the critical case.
\begin{corollary}\label{cor1} Let $d\ge 2$ and $p\in (2, 2^*)$. For $\,p\,$ sufficiently close to $\,2_+$, $\underline a(\vartheta(p,d),p)<a_\star^{\scriptscriptstyle\rm CKN}$, and so, there is symmetry breaking in a region where the radial extremals are linearly stable.
\end{corollary}
Notice that the case $d=2$ is not covered, for instance in Theorem~\ref{MainThm2} (ii), but the computations can be justified after noticing that among radial functions, \WLH\ also makes sense with $\gamma=d/2$ if $d=2$ (see \cite{DDFT}). By symmetry breaking, we mean $\C{CKN}^*(\vartheta(p,d),p,a)<\C{CKN}(\vartheta(p,d),p,a)$, since existence of extremals is not known for $a<a_\star^{\scriptscriptstyle\rm CKN}$.

\medskip\noindent {\it Proof of Proposition \ref{Lem:Limit}}. Optimal functions for Gagliardo-Nirenberg-Sobolev inequalities are, up to translations, radial solutions of the Euler-Lagrange equations
\be{EL-GN}
-\Delta u=\mathsf a\,u^{p-1}-\mathsf b\,u
\ee
where $\mathsf a$ and $\mathsf b$ are two positive coefficients which can be chosen arbitrarily because of the invariance of the inequality under a multiplication by a positive constant and the invariance under scalings. As a special choice, we can impose
\[
\mathsf a=\frac 2{p-2}\quad\mbox{and}\quad \mathsf b=\frac 2{p-2}-\frac d2\,(2+\log\pi)
\]
and denote by $u_p$ the (unique) corresponding solution so that, by passing to the limit as $p\to2_+$, we get the equation
\[
-\Delta u=2\,u\,\log u+\frac d2\,(2+\log\pi)\,u\,.
\]
Note that the function $u_2(x)$ is a positive radial solution of this equation in $\H^1(\R^d)$, which is normalized in $\L^2(\R^d)$: \hbox{$\nrm u2=1$}. According to \cite{MR1038450}, it is an extremal function for the logarithmic Sobolev inequality: for any $u\in\H^1(\R^d)
$,
\[
\ird{|u|^2\,\log\(\frac{|u|^2}{\nrm u2^2}\)}+\frac d2\,(2+\log\pi)\,\nrm u2^2\le\ird{|\nabla u|^2}\,,
\]
which, after optimization under scalings, is equivalent to~\eqref{LSI}. Moreover, it is unique as can be shown by considering for instance the remainder integral term arising from the Bakry-Emery method (see for instance \cite{MR1447044}, and \cite{MR1038450} for an earlier proof by a different method). A standard analysis shows that the solution $u_p$ converges to $u_2$ and \hbox{$\lim_{p\to2_+}\nrm{u_p}p=1$}. Multiplying~\eqref{EL-GN} by $u$ and by $x\cdot\nabla u$, one gets after a few integrations by parts that
\begin{eqnarray*}
&&\nrm{\nabla u_p}2^2=\mathsf a\,\nrm{u_p}p^p-\mathsf b\,\nrm{u_p}2^2\\
&&\frac{d-2}{2\,d}\,\nrm{\nabla u_p}2^2=\frac{\mathsf b}2\,\nrm{u_p}2^2-\frac{\mathsf a}p\,\nrm{u_p}p^p
\end{eqnarray*}
so that
\[
\textstyle\C{GN}(p)=\frac{\nrm{u_p}2^2}{\nrm{\nabla u_p}2^{2\,\vartheta(p,d)}\,\nrm{u_p}2^{2\,(1-\vartheta(p,d))}}=g(p)\,\nrm{u_p}p^{2-p}
\]
with $g(p):=\frac 12\(\frac{2\,p}d\)^{\vartheta(p,d)}\left[p\,\frac{4-d\,(p-2)\,(2+\log\pi)}{2\,p-d\,(p-2)}\right]^{1-\vartheta(p,d)}$ and the conclusion holds since $g'(2)=\frac d4\,\log\C{LS}$.
{\unskip\null\hfill$\square$\vskip 0.3cm}

Notice that it is possible to rephrase the Gagliardo-Nirenberg-Sobolev inequalities in a non scale invariant form as
\[
\mathsf a\,\nrm{u_p}p^{p-2}\,\nrm up^2\le\nrm{\nabla u}2^2+\mathsf b\,\nrm u2^2\quad\forall\;u\in\H^1(\R^d)\,,
\]
which can itself be recast into
\[
2\,\frac{\nrm up^2-\nrm u2^2}{p-2}\le\frac{\nrm{\nabla u}2^2}{\nrm{u_p}p^{p-2}}+\(\nrm{u_p}p^{2-p}\,\mathsf b-\mathsf a\)\nrm u2^2\,.
\]
It is then straightforward to understand why the limit as $p\to2_+$ in the above inequality gives the logarithmic Sobolev inequality (with optimal constant).

\bigskip\noindent{\scriptsize{\it Acknowledgements}.
The authors have been supported by the ANR projects CBDif-Fr and EVOL.\\[-4pt]\copyright\,2010 by the authors.\\[-4pt]
This paper may be reproduced, in its entirety, for non-commercial purposes.}\vspace*{-6pt}

\begin{thebibliography}{10}

\bibitem{Caffarelli-Kohn-Nirenberg-84}
{\sc L.~Caffarelli, R.~Kohn, and L.~Nirenberg}, {\em First order interpolation
  inequalities with weights}, Compositio Math., 53 (1984), pp.~259--275.

\bibitem{MR1038450}
{\sc E.~A. Carlen and M.~Loss}, {\em Extremals of functionals with competing
  symmetries}, J. Funct. Anal., 88 (1990), pp.~437--456.

\bibitem{Catrina-Wang-01}
{\sc F.~Catrina and Z.-Q. Wang}, {\em On the {C}affarelli-{K}ohn-{N}irenberg
  inequalities: sharp constants, existence (and nonexistence), and symmetry of
  extremal functions}, Comm. Pure Appl. Math., 54 (2001), pp.~229--258.

\bibitem{DDFT}
{\sc M.~del Pino, J.~Dolbeault, S.~Filippas, and A.~Tertikas}, {\em A
  logarithmic {H}ardy inequality}, Journal of Functional Analysis, 259 (2010),
  pp.~2045 -- 2072.

\bibitem{DE2010}
{\sc J.~Dolbeault and M.~J. Esteban}, {\em Extremal functions for
  {C}affarelli-{K}ohn-{N}irenberg and logarithmic {H}ardy inequalities}.
\newblock Preprint, 2010.

\bibitem{0902}
{\sc J.~Dolbeault, M.~J. Esteban, M.~Loss, and G.~Tarantello}, {\em On the
  symmetry of extremals for the {C}affarelli-{K}ohn-{N}irenberg inequalities},
  Adv. Nonlinear Stud., 9 (2009), pp.~713--726.

\bibitem{DET}
{\sc J.~Dolbeault, M.~J. Esteban, and G.~Tarantello}, {\em The role of {O}nofri
  type inequalities in the symmetry properties of extremals for
  {C}affarelli-{K}ohn-{N}irenberg inequalities, in two space dimensions}, Ann.
  Sc. Norm. Super. Pisa Cl. Sci.~(5), 7 (2008), pp.~313--341.

\bibitem{1007}
{\sc J.~Dolbeault, M.~J. Esteban, G.~Tarantello, and A.~Tertikas}, {\em Radial
  symmetry and symmetry breaking for some interpolation inequalities}.
\newblock Preprint, 2010.

\bibitem{Felli-Schneider-03}
{\sc V.~Felli and M.~Schneider}, {\em Perturbation results of critical elliptic
  equations of {C}affarelli-{K}ohn-{N}irenberg type}, J. Differential
  Equations, 191 (2003), pp.~121--142.

\bibitem{MR1447044}
{\sc G.~Toscani}, {\em Sur l'in\'egalit\'e logarithmique de {S}obolev}, C. R.
  Acad. Sci. Paris S\'er. I Math., 324 (1997), pp.~689--694.

\bibitem{MR479373}
{\sc F.~B. Weissler}, {\em Logarithmic {S}obolev inequalities for the
  heat-diffusion semigroup}, Trans. Amer. Math. Soc., 237 (1978), pp.~255--269.

\end{thebibliography}

\end{document}